# COLLATZ CONJECTURE: IS IT FALSE?

JUAN A. PEREZ

ABSTRACT. For a long time, Collatz Conjecture has been assumed to be true, although a formal proof has eluded all efforts to date. In this article, evidence is presented that suggests such an assumption is incorrect. By analysing the stopping times of various Collatz sequences, a pattern emerges that indicates the existence of non-empty sets of integers with stopping times greater than any given integer. This implies the existence of an infinite set of integers with non-finite stopping times, thus indicating the conjecture is false. Furthermore, a simple algorithm is constructed that finds integers with ever-greater stopping times. Such an algorithm does not halt, further supporting the conclusion that the conjecture is false.

## 1. INTRODUCTION

*Collatz Conjecture*, also known as *the $3x+1$ problem*, *the Syracuse problem*, *Kakutani's problem*, *Hasse's algorithm*, and *Ulam's problem* [2,3], concerns a simple arithmetic procedure applied to integers: If an integer $n$ is odd then "multiply by three and add one", while if it is even then "divide by two". Such an operation is described by the *Collatz function*

$$(1.1) \qquad C(n) = \begin{cases} 3n+1 & \text{if } n \equiv 1 \pmod{2}, \\ n/2 & \text{if } n \equiv 0 \pmod{2}. \end{cases}$$

The conjecture concerns the behaviour of this function under iteration, starting with any given positive integer $n$: it is stated that, in all cases, successive iterations of the function $C(n)$ will eventually reach the number 1; thereafter iterations will cycle, taking successive values 1, 4, 2, 1, ...

Although arithmetically simple to describe, Collatz Conjecture remains an unsolved problem which appears to be extremely difficult to prove [2,3]. As such, it has fascinated mathematicians and non-mathematicians alike. Over the years, it has been studied by mathematicians, physicists and computer scientists, generating hundreds of publications [3,6].

This article presents evidence strongly suggesting the conjecture is false. By reformulating the Collatz function $C(n)$, it focuses on the sequences of odd numbers generated by the iterations and the so-called *stopping time*, $\sigma(n)$, i.e. the number of iterations (steps) required for $C(n) < n$. These stopping times build a pattern that implies the existence of non-empty sets of positive integers ($n$) with stopping times greater than any given integer $k$, $\sigma(n) > k$. This in turn implies the existence of an infinite set of integers with non-finite stopping times.

## 2. COLLATZ SEQUENCES

When investigating the behaviour of the Collatz function $C(n)$, it is of value to consider the factors that determine the length of the corresponding *Collatz sequences*, i.e. the sequences of terms $C_i(n)$ generated by the $i$ iterations required to reach firstly their corresponding *stopping time* $\sigma(n)$, and secondly their subsequent *total stopping time* $\sigma_\infty(n)$, i.e. the number of iterations necessary to reach 1.

It is trivial to observe that, when $n$ is an even integer, $\sigma(n)$ is always 1. Consequently, when studying $\sigma(n)$, it becomes more convenient to consider only odd integers and their corresponding Collatz sequences. Originally used by Crandall [1], a different function $O(n)$ is proposed, such that

$$(2.1) \qquad O(n) := (3n+1)/2^m \equiv 1 \,(\mathrm{mod}\, 2)$$

where $2^m$ is the power of 2 that makes $O(n)$ an odd integer. As an example, the Collatz sequence for the integer $n = 15$ will be

$$15 \xrightarrow[/2]{\times 3+1} 23 \xrightarrow[/2]{\times 3+1} 35 \xrightarrow[/2]{\times 3+1} 53 \xrightarrow[/2^5]{\times 3+1} 5 \xrightarrow[/2^4]{\times 3+1} 1$$

hence making $\sigma(15) = 4$ and $\sigma_\infty(15) = 5$. The purpose of condensing the two operations ($\times 3+1$ and $/2^m$) into one single step is to generate Collatz sequences with exclusively odd numbers.

It is a very simple observation that, for the Collatz Conjecture to be true, it is both necessary and sufficient that all odd positive integers have finite stopping times. This will invariably make all Collatz sequences reach 1. Therefore, an evaluation of the factors affecting the behaviour of the sequences with regard to their stopping times will shed light on the validity of the conjecture.

The first easy observation to be made is that, for all odd integers $n \equiv 1 \,(\mathrm{mod}\, 4)$, stopping times are always $\sigma(n) = 1$. Since $3n + 1 \equiv 0 \,(\mathrm{mod}\, 4)$, the first term in the sequence, $O_1(n)$, will be derived by dividing by at least $2^2$. Therefore,

TABLE 1. Patterns of stopping times, σ(n), for an initial set of odd integers, n.

| Step | n |
|---|---|
| | 1 3 5 7 9 11 13 15 17 19 21 23 25 27 29 31 33 35 37 39 41 43 45 47 |
| 1 | + − + − + − + − + − + − + − + − + − + − + − + − |

| Step | n |
|---|---|
| | 3 7 11 15 19 23 27 31 35 39 43 47 51 55 59 63 67 71 75 79 83 87 91 95 |
| 2 | + − − − + − − − + − − − + − − − + − − − + − − − |
| 3 | + − + − + + − − + − + − + + − − + − + − + + − − |

"+": Sequence has reached its stopping time. "−": Sequence still to reach its stopping time.

it can be written that $O_1(n) = (3n+1)/4 = 0.75n + 0.25$. Since it can also be written that $n = 2p + 1$, this makes $O_1(n) = 0.75(2p+1) + 0.25 = 1.5p + 1$, hence $O_1(n) < n$.

Alternatively, when the odd integers are $n \equiv 3 \pmod{4}$, stopping times are always greater than 1. Since $3n + 1 \equiv 2 \pmod{4}$, the term $3n + 1$ can only be divided by a single 2. Hence, $O_1(n) = (3n+1)/2 = 1.5n + 0.5$, and $O_1(n) > n$.

The above observations generate an initial pattern in the stopping times for all odd integers, alternating between $\sigma(n) = 1$ and $\sigma(n) > 1$: for the first step in the Collatz sequences, 50% of all odd integers reach their stopping times while the other half do not. This suggests that it will be of interest to look into any possible pattern emerging as the sequences progress towards more and more steps. Table 1 shows what happens for the first three steps in the Collatz sequences. The criterion used is to label "+" when the stopping time has been reached, and "−" when it has not. For the first step, the pattern is "+ −", as already explained. For the second step, the table shows the pattern observed for odd integers $n \equiv 3 \pmod{4}$, "+ − − −". Such a pattern is repeated along the number line. For the third step, the pattern found is "+ − + − + + − −", which is also repeated along the number line.

For step 1, the pattern indicates that the stopping times in the whole number line obey the template found for odd integers $n \equiv 1, 3 \pmod{2^2}$. For step 2, the same happens for the template of odd integers $n \equiv 1, 3, 5, 7, 9, 11, 13, 15 \pmod{2^4}$. And for step 3, the corresponding template is set for odd congruences $\pmod{2^5}$. For subsequent steps, similar patterns emerge, with templates set by the odd congruences mod $2^7$, mod $2^8$, mod $2^{10}$, mod $2^{12}$, and so on.

To understand the reason behind these patterns, it is necessary to examine the impact that the factors $2^m$ have on the values of the terms $O(n)$ in (2.1). For this purpose, consider the Collatz sequence of a given example, $n = 191$, as detailed in Table 2. As well as the values of $O_i(n)$ for every step $i$ up to the stopping time $\sigma(n) = 8$, the table also includes the values of $m_i$ for every step (third column in the table). The fourth column shows the "accumulated 2s" ($\Sigma m_i$) down the sequence. These values are compared to the "maximal 2s" (fifth column in the table), i.e. the total number of 2s that would have been required for the sequence to reach its stopping time at that point. Since every term in the sequence is generated multiplying by 3 (and adding 1) and dividing by $2^{m_i}$, the stopping time will be reached when $3^i < 2^{\Sigma m_i}$. In the example used in the table, this happens after 8 steps, when the accumulated total of 2s reaches (or exceeds) the maximal number of 2s required for $3^i < 2^{\Sigma m_i}$. In general, if $\sigma(n) = x$, the required total of accumulated 2s, $\Sigma m_i = y$, will be determined by $3^x < 2^y$, this is

$$(2.2) \qquad y = \left\lceil \frac{x \operatorname{Ln} 3}{\operatorname{Ln} 2} \right\rceil$$

TABLE 2. Collatz sequence of integer $n=191$ up to its stopping time $\sigma(n)=8$.

| Step (i) | Integer $O_i(n)$ | 2s in step ($m_i$) | Accumulated 2s ($\Sigma m_i$) | Maximal 2s | 2s deficit |
|---|---|---|---|---|---|
| 0 | 191 | – | – | – | – |
| 1 | 287 | 1 | 1 | 2 | 1 |
| 2 | 431 | 1 | 2 | 4 | 2 |
| 3 | 647 | 1 | 3 | 5 | 2 |
| 4 | 971 | 1 | 4 | 7 | 3 |
| 5 | 1,457 | 1 | 5 | 8 | 3 |
| 6 | 1,093 | 2 | 7 | 10 | 3 |
| 7 | 205 | 4 | 11 | 12 | 1 |
| 8 | 77 | 3 | 14 | 13 | -1 |

In the example used in Table 2, $x=8$ and $y=13$. The sixth column in the table shows the "2s deficit", i.e. the difference between the maximal total of 2s required to reach the stopping time and the accumulated 2s. Only when this deficit becomes 0 (or negative) is the stopping time reached.

The dynamics of the Collatz sequence highlighted in Table 2 are universal and can be used to explain the patterns described previously (Table 1). Two simple theorems are sufficient.

**Theorem 2.1.** *Consider the Collatz sequence of an integer $n_0$ with stopping time $\sigma(n_0) = x$. Consider also a new integer $n = n_0 + k 2^y$, where k can be any positive integer and y is given by* (2.2). *The stopping time of the new integer will always be equal to that of $n_0$, $\sigma(n) = \sigma(n_0) = x$.*

*Proof.* The terms of the Collatz sequence for $n$, $O_i(n)$, are as follows:

$$O_1(n) = O_1(n_0) + 3k2^{y-m_1}$$
$$O_2(n) = O_2(n_0) + 3^2 k 2^{y-(m_1+m_2)}$$
$$\ldots\ldots$$
$$O_x(n) = O_x(n_0) + 3^x k$$

Given that $O_x(n_0) < n_0$ and $k3^x < k2^y$, it can be concluded that $O_x(n) < n$ and, consequently, $\sigma(n) = \sigma(n_0) = x$. □

In the proof of Theorem 2.1, it should be noted that all the intermediate terms $3^i k 2^{y-\Sigma m_i}$ are greater than $k2^y$, since not enough 2s have been eliminated to reach the stopping time.

**Theorem 2.2.** *Consider the Collatz sequence of an integer $n_0$ with stopping time $\sigma(n_0) > x$. Consider also a new integer $n = n_0 + k2^y$, where $k$ can be any positive integer and $y$ is given by* (2.2). *The stopping time of the new integer will also be greater than $x$, $\sigma(n) > x$.*

*Proof.* The terms of the Collatz sequence for $n$, $O_i(n)$, are as follows:

$$O_1(n) = O_1(n_0) + 3k2^{y-m_1}$$
$$O_2(n) = O_2(n_0) + 3^2 k 2^{y-(m_1+m_2)}$$
$$\ldots\ldots$$
$$O_x(n) = O_x(n_0) + 3^x k 2^{y-\Sigma m_i}$$

Given that $O_x(n_0) > n_0$ and $3^x k 2^{y-\Sigma m_i} > k2^y$, it is concluded that $O_x(n) > n$ and, consequently, $\sigma(n) > x$. □

Theorems 2.1 and 2.2 are all that is required to explain the patterns found for the stopping times. This means that the templates (mod $2^y$) observed in Table 1 repeat along the number line, not only for steps 1, 2 and 3, but for all subsequent steps. The values of $y$, as given by (2.2), are directly determined by the stopping times, $\sigma(n) = x$ (see Table 3). The patterns generated are a fundamental property of the stopping times and have strong implications for Collatz Conjecture.

## 3. CONVERSION RATES

A way of looking into the validity of Collatz Conjecture is to analyse the conversion rates of the Collatz sequences, i.e. the rates at which the stopping times are reached. It has already been proved by others that "almost every" positive integer has a finite stopping time [3,5]. This suggests that the Collatz sequences will convert from "–" to "+" (following the convention used in Table 1) at a rate high enough to guarantee the eventual and full conversion of all the sequences, so they all have finite stopping times. The periodicity of the patterns described in the previous section permits calculation of the global rate of conversion from one step to the next. Table 4 shows the results obtained from such a calculation.

TABLE 3. Length of templates (mod $2^y$).

| σ(n)=x | x in $3^x$ | y in $2^y$ | $2^y$ |
|---|---|---|---|
| 1 | 1 | 2 | 4 |
| 2 | 2 | 4 | 16 |
| 3 | 3 | 5 | 32 |
| 4 | 4 | 7 | 128 |
| 5 | 5 | 8 | 256 |
| 6 | 6 | 10 | 1,024 |
| 7 | 7 | 12 | 4,096 |
| 8 | 8 | 13 | 8,192 |
| 9 | 9 | 15 | 32,768 |
| 10 | 10 | 16 | 65,536 |
| 11 | 11 | 18 | 262,144 |
| 12 | 12 | 20 | 1,048,576 |
| 13 | 13 | 21 | 2,097,152 |

The left-hand side of Table 4 shows how the percentage of sequences still to reach the stopping time decreases gradually as the sequences progress from one step to the next. This is, in principle, compatible with the proven fact that "almost every" positive integer has a finite stopping time [3,5], since the density of sequences with unreached stopping times decreases apparently towards zero. Figure 1 shows the same results in graphical form.

TABLE 4. Conversion rates for the initial steps of the Collatz sequences.

| Step | Unreached stopping time (out of total odd numbers) | % | Unreached stopping time (out of remaining non-stopped sequences) | Non-conversion rate (%) |
|---|---|---|---|---|
| 1 | 1 (2) | 50.0 | 1 (2) | 50.0 |
| 2 | 3 (8) | 37.5 | 3 (4) | 75.0 |
| 3 | 4 (16) | 25.0 | 4 (6) | 66.7 |
| 4 | 13 (64) | 20.3 | 13 (16) | 81.2 |
| 5 | 19 (128) | 14.8 | 19 (26) | 73.1 |
| 6 | 64 (512) | 12.5 | 64 (76) | 84.2 |
| 7 | 226 (2,048) | 11.0 | 226 (256) | 88.3 |
| 8 | 367 (4,096) | 9.0 | 367 (452) | 81.2 |
| 9 | 1,294 (16,384) | 7.9 | 1,294 (1,468) | 88.1 |
| 10 | 2,114 (32,768) | 6.5 | 2,114 (2,588) | 81.7 |

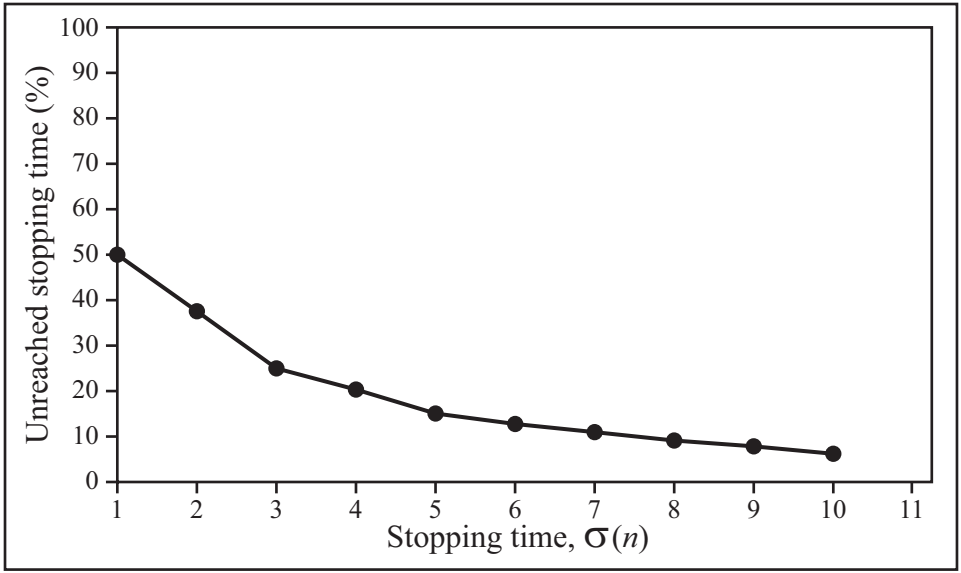

FIGURE 1. Density of sequences with unreached stopping times.

However, the right-hand side of Table 4 tells a different story. The repeated patterns for each step in the the Collatz sequences permit calculation of the rate of conversion from unreached to reached stopping times. These rates are surprisingly low. As a consequence, the rates of non-conversion (i.e. the percentage of sequences still to reach their stopping times) take values in excess of 80% (Figure 2 shows the same results in graphical form). This means that the total number of sequences with unreached stopping times for each template (mod $2^y$) increases as the sequences progress from one step to the next. The jump from one template to the next (as seen in Table 3) is always a factor of 2 or 4, so the conversion rates observed are insufficient to establish a diminishing trend for the unreached stopping times.

The implications of these results for Collatz Conjecture are considerable. They clearly establish the existence, for every step in the sequences (i.e. for every positive integer $x$), of an infinite set of integers $n$ with stopping times $\sigma(n) > x$. Assuming that the non-conversion rates recorded in Table 4 and Figure 2 remain at the same level for ever greater values of $x$, the eventual consequence will be the existence of an infinite set of integers with stopping times greater than any value of $x$, i.e. non-finite. In other words, Collatz Conjecture will be false.

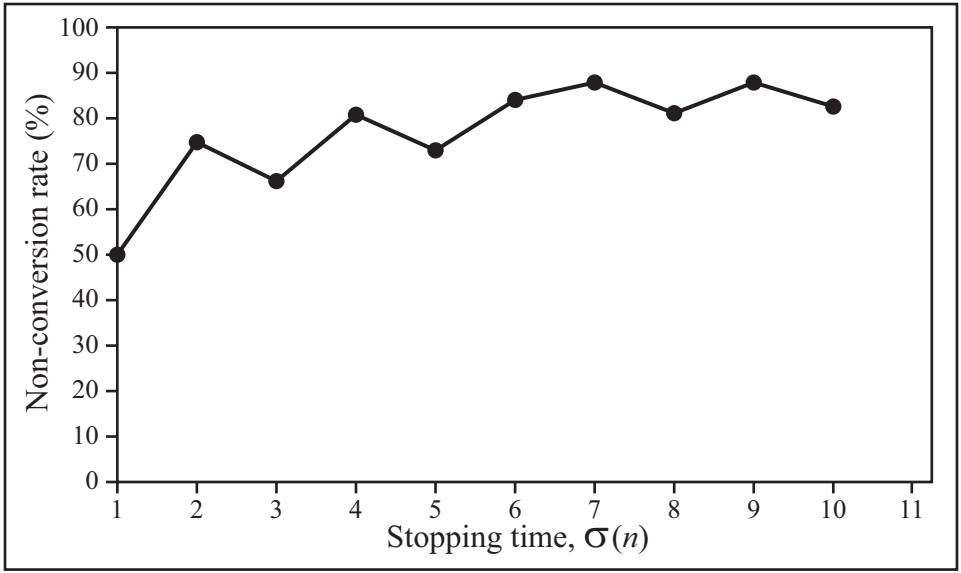

FIGURE 2. Non-conversion rates of unreached stopping times.

It can be argued that non-conversion rates greater than 50% are to be expected. A number of factors have to be taken into consideration:
- The final step of a Collatz sequence that reaches its stopping time will always involve a power of 2 of at least 2 units, i.e. $2^{m_i}=2^2$, or greater.
- The "2s deficit" (as introduced in Table 2), carried over by the sequence, will have to be sufficiently small for the last factor $2^{m_i}$ to reach (or surpass) the maximal value $2^y$ set by (2.2).
- For any step of a Collatz sequence, the value of the corresponding factor $2^{m_i}$ will be, on average, $2^1$ (50% of cases), $2^2$ (25%), $2^3$ (12.5%), $2^4$ (6.3%), $2^5$ (3.1%), and so on. This is determined by the prime number factorization of even numbers: 50% of all even numbers are divisible only by $2^1$, 25% are divisible by $2^2$, 12.5% are divisible by $2^3$, and so on.
- Consequently, when considering the set of Collatz sequences still to reach the stopping time after a given number of steps, around 50% of the total will face a division by only $2^1$, hence they will remain unable to reach it.
- The remaining 50% of sequences will face a division by a higher power of 2, i.e. 25% by $2^2$, 12.5% by $2^3$, and so on. But they will only be able to reach their stopping time if their "2s deficit" is sufficiently small – the fact that the non-conversion rates observed are greater than 80% indicates that less than 20% of sequences satisfy this requirement.

These considerations assume that the occurrence of even numbers in the Collatz sequences with the divisibilities by powers of 2 described here is similar to that of the whole set of positive integers. Table 5 shows the observed occurrence of powers of 2 ($2^{m_i}$) in some selected Collatz sequences, up to their stopping time $\sigma(n)$, or their total stopping time $\sigma_\infty(n)$. It can be seen that they are in broad agreement with the theoretical percentages, confirming the validity of the assumption.

TABLE 5. Occurrence of powers of 2 ($2^{m_i}$) in selected Collatz sequences.

| Integer ($n$) | Stopping time, $\sigma(n)$ | $2^1$s (%) | $2^2$s (%) | $2^3$s (%) | $2^4$s (%) | $2^5$s (%) | $2^6$s (%) |
|---|---|---|---|---|---|---|---|
| 27 | 37 | 22 (59.5) | 10 (27.0) | 3 (8.1) | 2 (5.4) | 0 (0.0) | 0 (0.0) |
| 1055 | 50 | 32 (64.0) | 11 (22.0) | 3 (6.0) | 2 (4.0) | 1 (2.0) | 1 (2.0) |
| 7279 | 48 | 30 (62.5) | 12 (25.0) | 3 (6.3) | 0 (0.0) | 2 (4.2) | 1 (2.1) |
| Theoretical percentages | | 50.0 | 25.0 | 12.5 | 6.3 | 3.1 | 1.6 |

| Integer ($n$) | Total stopping time, $\sigma_\infty(n)$ | $2^1$s (%) | $2^2$s (%) | $2^3$s (%) | $2^4$s (%) | $2^5$s (%) | $2^6$s (%) | $2^7$s (%) |
|---|---|---|---|---|---|---|---|---|
| 27 | 41 | 24 (58.5) | 10 (24.4) | 3 (7.3) | 3 (7.3) | 1 (2.4) | 0 (0.0) | 0 (0.0) |
| $27+2^{51}$ | 191 | 104 (54.5) | 50 (26.2) | 18 (9.4) | 8 (4.2) | 6 (3.1) | 3 (1.6) | 1 (0.5) |
| Theoretical percentages | | 50.0 | 25.0 | 12.5 | 6.3 | 3.1 | 1.6 | 0.8 |

There are strong reasons to conclude that the non-conversion rates recorded in Table 4 and Figure 2 remain at the same level as the number of steps in the Collatz sequences are increased indefinitely. Hence, it can be deduced that the set of positive integers $n$ with stopping times greater than any given integer $x$ (i.e. non-finite) is infinite in size. This implies Collatz Conjecture is false.

## 4. LOOKING FOR COUNTER-EXAMPLES

Although the evidence presented in this article refutes Collatz Conjecture, it is also compatible with previous results regarding the behaviour of the Collatz sequences and their stopping times [3,5]. It has already been reported that "almost every" positive integer has a finite stopping time [5]. The results presented here are compatible with that statement: as shown in Table 4 and Figure 1, the total occurrence of sequences with unreached stopping times (as a percentage of all odd integers up to the corresponding modulo $2^y$) decreases gradually as the number of steps increases. Given that modulo $2^y$ grows exponentially, such a percentage will be expected to reach very low values, approaching zero.

But the non-conversion rates recorded in Table 4 and Figure 2 also imply that the number of integers with stopping times greater than a given integer $x$ actually increases from one modulo $2^y$ to the next. As already pointed out, the consequence of such a trend is the existence of an infinite set of integers $n$ with stopping times greater than *any* integer $x$ (i.e. non-finite), hence implying Collatz Conjecture is false.

If the conjecture is false, is it possible to find integers $n$ with a non-finite stopping time, $\sigma(n) = \infty$? Although the existence of such integers can be deduced, the answer is probably negative. Not only are such counter-examples likely to be extremely large (the conjecture has already been found to hold true for all $n < 20 \times 2^{58}$ [5]), proving that a given integer $n$ has a non-finite stopping time will require to show conclusively that the corresponding Collatz sequence *never* reaches a point where $\sigma(n) < n$, something that might not be at all feasible.

Nevertheless, an algorithm was constructed to find integers $n$ with ever-increasing stopping times. The following set of instructions was implemented:

*a*) Select a small integer $n_0$ and examine its Collatz sequence (as in Table 2), looking for a step as close as possible to the stopping time where the jump from a modulo $2^y$ to the next modulo involves a factor of 4.

*b*) Construct a new integer $n_1$ by taking $n_0$ and adding to it $2^y$, $2 \times 2^y$, or $3 \times 2^y$.

*c*) Inspect the Collatz sequences for the three newly constructed integers (i.e. $n_0 + 2^y$, $n_0 + 2 \times 2^y$, $n_0 + 3 \times 2^y$), and select the sequence with the highest value for the stopping time $\sigma(n_1)$ of the three.

*d*) Inspect the Collatz sequence for the selected integer $n_1$ and look again for a step as close as possible to the stopping time where the jump between modulos $2^y$ is a factor of 4.

*e*) Construct a new integer $n_2$ by taking $n_1$ and adding to it the new values $2^y$, $2 \times 2^y$, or $3 \times 2^y$.

*f*) Of the three Collatz sequences, select again the sequence with the greatest stopping time $\sigma(n_2)$.

*g*) Repeat the process of constructing new integers and selecting the sequence with the greatest stopping time.

*h*) If at any point in the iterations, none of the three constructed integers $n_i$ produces a stopping time greater than that of the preceding integer $n_{i-1}$, try constructing a new (three new) integer(s) by selecting the modulo $2^y$ (step)

TABLE 6. Ever-increasing stopping times, $\sigma(n_i)$, for a set of integers, $n_i$.

| $n_i$ | $\sigma(n_i)$ | $n_i$ | $\sigma(n_i)$ | $n_i$ | $\sigma(n_i)$ |
|---|---|---|---|---|---|
| $n_0 = 27$ | 37 | $n_{54} = n_{53} + 2 \times 2^{934}$ | 598 | $n_{108} = n_{107} + 3 \times 2^{1,918}$ | 1,218 |
| $n_1 = n_0 + 2^{57}$ | 48 | $n_{55} = n_{54} + 2 \times 2^{945}$ | 600 | $n_{109} = n_{108} + 2^{1,929}$ | 1,222 |
| $n_2 = n_1 + 3 \times 2^{75}$ | 51 | $n_{56} = n_{55} + 3 \times 2^{948}$ | 601 | $n_{110} = n_{109} + 3 \times 2^{1,934}$ | 1,229 |
| $n_3 = n_2 + 2^{78}$ | 52 | $n_{57} = n_{56} + 2^{951}$ | 619 | $n_{111} = n_{110} + 2 \times 2^{1,945}$ | 1,242 |
| $n_4 = n_3 + 3 \times 2^{81}$ | 59 | $n_{58} = n_{57} + 2^{980}$ | 624 | $n_{112} = n_{111} + 2 \times 2^{1,967}$ | 1,249 |
| $n_5 = n_4 + 3 \times 2^{92}$ | 92 | $n_{59} = n_{58} + 3 \times 2^{988}$ | 649 | $n_{113} = n_{112} + 2 \times 2^{1,977}$ | 1,281 |
| $n_6 = n_5 + 3 \times 2^{143}$ | 101 | $n_{60} = n_{59} + 2^{1,026}$ | 660 | $n_{114} = n_{113} + 2^{2,029}$ | 1,285 |
| $n_7 = n_6 + 2 \times 2^{158}$ | 107 | $n_{61} = n_{60} + 3 \times 2^{1,045}$ | 680 | $n_{115} = n_{114} + 2^{2,034}$ | 1,286 |
| $n_8 = n_7 + 2 \times 2^{167}$ | 112 | $n_{62} = n_{61} + 3 \times 2^{1,075}$ | 681 | $n_{116} = n_{115} + 3 \times 2^{2,037}$ | 1,295 |
| $n_9 = n_8 + 2 \times 2^{176}$ | 119 | $n_{63} = n_{62} + 2 \times 2^{1,078}$ | 689 | $n_{117} = n_{116} + 2^{2,051}$ | 1,302 |
| $n_{10} = n_9 + 2^{186}$ | 120 | $n_{64} = n_{63} + 2^{1,091}$ | 696 | $n_{118} = n_{117} + 2^{2,061}$ | 1,307 |
| $n_{11} = n_{10} + 3 \times 2^{189}$ | 125 | $n_{65} = n_{64} + 2^{1,102}$ | 710 | $n_{119} = n_{118} + 2^{2,070}$ | 1,313 |
| $n_{12} = n_{11} + 2^{197}$ | 138 | $n_{66} = n_{65} + 2^{1,124}$ | 723 | $n_{120} = n_{119} + 2^{2,080}$ | 1,326 |
| $n_{13} = n_{12} + 2^{216}$ | 141 | $n_{67} = n_{66} + 2^{1,143}$ | 740 | $n_{121} = n_{120} + 3 \times 2^{2,099}$ | 1,331 |
| $n_{14} = n_{13} + 2^{222}$ | 148 | $n_{68} = n_{67} + 2 \times 2^{1,170}$ | 743 | $n_{122} = n_{121} + 2^{2,107}$ | 1,345 |
| $n_{15} = n_{14} + 2 \times 2^{233}$ | 158 | $n_{69} = n_{68} + 2^{1,175}$ | 744 | $n_{123} = n_{122} + 2^{2,129}$ | 1,368 |
| $n_{16} = n_{15} + 2^{249}$ | 160 | $n_{70} = n_{69} + 3 \times 2^{1,178}$ | 751 | $n_{124} = n_{123} + 3 \times 2^{2,167}$ | 1,383 |
| $n_{17} = n_{16} + 3 \times 2^{251}$ | 163 | $n_{71} = n_{70} + 2 \times 2^{1,189}$ | 762 | $n_{125} = n_{124} + 2^{2,191}$ | 1,392 |
| $n_{18} = n_{17} + 2^{257}$ | 177 | $n_{72} = n_{71} + 3 \times 2^{1,205}$ | 773 | $n_{126} = n_{125} + 3 \times 2^{2,205}$ | 1,403 |
| $n_{19} = n_{18} + 2 \times 2^{279}$ | 183 | $n_{73} = n_{72} + 2^{1,224}$ | 777 | $n_{127} = n_{126} + 2^{2,221}$ | 1,407 |
| $n_{20} = n_{19} + 3 \times 2^{289}$ | 195 | $n_{74} = n_{73} + 2^{1,230}$ | 782 | $n_{128} = n_{127} + 2^{2,229}$ | 1,410 |
| $n_{21} = n_{20} + 3 \times 2^{308}$ | 198 | $n_{75} = n_{74} + 3 \times 2^{1,238}$ | 792 | $n_{129} = n_{128} + 3 \times 2^{2,232}$ | 1,436 |
| $n_{22} = n_{21} + 3 \times 2^{311}$ | 211 | $n_{76} = n_{75} + 3 \times 2^{1,254}$ | 800 | $n_{130} = n_{129} + 3 \times 2^{2,275}$ | 1,449 |
| $n_{23} = n_{22} + 2^{333}$ | 226 | $n_{77} = n_{76} + 3 \times 2^{1,265}$ | 803 | $n_{131} = n_{130} + 2 \times 2^{2,294}$ | 1,459 |
| $n_{24} = n_{23} + 2 \times 2^{357}$ | 244 | $n_{78} = n_{77} + 2 \times 2^{1,270}$ | 817 | $n_{132} = n_{131} + 2^{2,311}$ | 1,474 |
| $n_{25} = n_{24} + 3 \times 2^{384}$ | 256 | $n_{79} = n_{78} + 2^{1,292}$ | 823 | $n_{133} = n_{132} + 3 \times 2^{2,335}$ | 1,492 |
| $n_{26} = n_{25} + 3 \times 2^{403}$ | 258 | $n_{80} = n_{79} + 2^{1,303}$ | 828 | $n_{134} = n_{133} + 3 \times 2^{2,362}$ | 1,497 |
| $n_{27} = n_{26} + 3 \times 2^{406}$ | 264 | $n_{81} = n_{80} + 3 \times 2^{1,311}$ | 841 | $n_{135} = n_{134} + 2^{2,370}$ | 1,507 |
| $n_{28} = n_{27} + 3 \times 2^{417}$ | 270 | $n_{82} = n_{81} + 3 \times 2^{1,330}$ | 852 | $n_{136} = n_{135} + 3 \times 2^{2,387}$ | 1,530 |
| $n_{29} = n_{28} + 3 \times 2^{425}$ | 272 | $n_{83} = n_{82} + 3 \times 2^{1,349}$ | 870 | $n_{137} = n_{136} + 3 \times 2^{2,422}$ | 1,546 |
| $n_{30} = n_{29} + 3 \times 2^{430}$ | 280 | $n_{84} = n_{83} + 3 \times 2^{1,376}$ | 875 | $n_{138} = n_{137} + 3 \times 2^{2,449}$ | 1,569 |
| $n_{31} = n_{30} + 3 \times 2^{441}$ | 282 | $n_{85} = n_{84} + 2^{1,384}$ | 894 | $n_{139} = n_{138} + 3 \times 2^{2,484}$ | 1,594 |
| $n_{32} = n_{31} + 3 \times 2^{444}$ | 284 | $n_{86} = n_{85} + 2^{1,414}$ | 931 | $n_{140} = n_{139} + 3 \times 2^{2,525}$ | 1,622 |
| $n_{33} = n_{32} + 3 \times 2^{449}$ | 307 | $n_{87} = n_{86} + 2^{1,473}$ | 933 | $n_{141} = n_{140} + 3 \times 2^{2,568}$ | 1,625 |
| $n_{34} = n_{33} + 2^{485}$ | 358 | $n_{88} = n_{87} + 2^{1,476}$ | 942 | $n_{142} = n_{141} + 2^{2,574}$ | 1,630 |
| $n_{35} = n_{34} + 2 \times 2^{566}$ | 374 | $n_{89} = n_{88} + 3 \times 2^{1,492}$ | 952 | $n_{143} = n_{142} + 3 \times 2^{2,582}$ | 1,735 |
| $n_{36} = n_{35} + 3 \times 2^{590}$ | 383 | $n_{90} = n_{89} + 3 \times 2^{1,506}$ | 986 | $n_{144} = n_{143} + 2 \times 2^{2,747}$ | 1,738 |
| $n_{37} = n_{36} + 2^{606}$ | 394 | $n_{91} = n_{90} + 2 \times 2^{1,560}$ | 988 | $n_{145} = n_{144} + 2^{2,752}$ | 1,745 |
| $n_{38} = n_{37} + 3 \times 2^{623}$ | 438 | $n_{92} = n_{91} + 3 \times 2^{1,563}$ | 991 | $n_{146} = n_{145} + 2 \times 2^{2,763}$ | 1,764 |
| $n_{39} = n_{38} + 2^{693}$ | 441 | $n_{93} = n_{92} + 3 \times 2^{1,568}$ | 1,013 | $n_{147} = n_{146} + 3 \times 2^{2,793}$ | 1,785 |
| $n_{40} = n_{39} + 2^{696}$ | 448 | $n_{94} = n_{93} + 2^{1,604}$ | 1,022 | $n_{148} = n_{147} + 2^{2,828}$ | 1,788 |
| $n_{41} = n_{40} + 2^{709}$ | 455 | $n_{95} = n_{94} + 3 \times 2^{1,617}$ | 1,023 | $n_{149} = n_{148} + 2^{2,831}$ | 1,814 |
| $n_{42} = n_{41} + 2^{720}$ | 463 | $n_{96} = n_{95} + 2^{1,620}$ | 1,036 | $n_{150} = n_{149} + 2^{2,874}$ | 1,816 |
| $n_{43} = n_{42} + 2^{731}$ | 464 | $n_{97} = n_{96} + 2^{1,641}$ | 1,041 | $n_{151} = n_{150} + 3 \times 2^{2,877}$ | 1,842 |
| $n_{44} = n_{43} + 2^{734}$ | 467 | $n_{98} = n_{97} + 3 \times 2^{1,647}$ | 1,055 | $n_{152} = n_{151} + 3 \times 2^{2,918}$ | 1,858 |
| $n_{45} = n_{44} + 2^{739}$ | 470 | $n_{99} = n_{98} + 3 \times 2^{1,671}$ | 1,075 | $n_{153} = n_{152} + 3 \times 2^{2,942}$ | 1,865 |
| $n_{46} = n_{45} + 3 \times 2^{742}$ | 475 | $n_{100} = n_{99} + 2 \times 2^{1,701}$ | 1,139 | $n_{154} = n_{153} + 2 \times 2^{2,953}$ | 1,876 |
| $n_{47} = n_{46} + 2^{750}$ | 483 | $n_{101} = n_{100} + 3 \times 2^{1,804}$ | 1,148 | $n_{155} = n_{154} + 3 \times 2^{2,972}$ | 1,890 |
| $n_{48} = n_{47} + 3 \times 2^{764}$ | 495 | $n_{102} = n_{101} + 3 \times 2^{1,818}$ | 1,155 | $n_{156} = n_{155} + 3 \times 2^{2,994}$ | 1,904 |
| $n_{49} = n_{48} + 2^{783}$ | 534 | $n_{103} = n_{102} + 2^{1,828}$ | 1,163 | $n_{157} = n_{156} + 2 \times 2^{3,015}$ | 1,919 |
| $n_{50} = n_{49} + 2 \times 2^{845}$ | 551 | $n_{104} = n_{103} + 2^{1,842}$ | 1,178 | $n_{158} = n_{157} + 2^{3,040}$ | 1,969 |
| $n_{51} = n_{50} + 3 \times 2^{872}$ | 575 | $n_{105} = n_{104} + 2^{1,866}$ | 1,180 | $n_{159} = n_{158} + 2 \times 2^{3,118}$ | 2,003 |
| $n_{52} = n_{51} + 2 \times 2^{910}$ | 579 | $n_{106} = n_{105} + 3 \times 2^{1,869}$ | 1,204 | $n_{160} = n_{159} + 2 \times 2^{3,172}$ | 2,012 |
| $n_{53} = n_{52} + 2 \times 2^{915}$ | 591 | $n_{107} = n_{106} + 2^{1,907}$ | 1,211 | | |

that immediately precedes (or follows) the step used previously. As before, select the sequence with the greatest stopping time $\sigma(n_i)$ and continue the iterations.

Table 6 shows the results obtained by following the algorithm just described, starting with $n_0 = 27$. Of the 160 iterations used to construct the table, on no occasion was it necessary to look for alternatives to the three newly constructed integers $n_i + 2^y$, $n_i + 2 \times 2^y$ and $n_i + 3 \times 2^y$.

The rationale of the algorithm is based on the patterns and rates of non-conversion reported in this article. When the stopping time $\sigma(n_i)$ is reached in a given Collatz sequence, the preceding step holds the key to greater stopping times. If the jump between both steps is a factor of 4, there will be four sequences to consider, corresponding to the integers $n_i$, $n_i + 2^y$, $n_i + 2 \times 2^y$ and $n_i + 3 \times 2^y$. The sequence for $n_i$ will have $\sigma(n_i)$ as stopping time. But the remaining three sequences will obey the rates of non-conversion reported in this article, i.e. in excess of 80%. Therefore the probability of at least one of the sequences having a stopping time greater than $\sigma(n_i)$ will be very high. Table 6 validates this assessment.

Table 6 also demonstrates the reliability of the algorithm devised to search for greater stopping times. It strongly supports the claim already made in this article that the high rates of non-conversion recorded in Table 4 and Figure 2 are maintained as progress is made towards greater stopping times. If this had not been the case, the algorithm would not have been able to find integers $n_i$ with ever greater stopping times $\sigma(n_i)$ as readily as observed. The efficiency of the algorithm does not diminish as it progresses towards more and more iterations – this can be seen in Figure 3, where the stopping times of Table 6 are plotted against the number of iterations. The implication is that, if left to run indefinitely, the algorithm will never halt.

It must be observed that it is in the nature of this algorithm to use larger and larger integers $n_i$ in the search for greater stopping times $\sigma(n_i)$. A linear relationship between $\sigma(n_i)$ and $\text{Log}_2(n_i)$ is also found (see Figure 4), which can be explained as follows. By rearranging (2.2), it can be written that

$$y \, \text{Ln} 2 \simeq x \, \text{Ln} 3$$

where $x = \sigma(n_i)$ and $y \simeq \text{Log}_2(n_i)$. Therefore,

(4.1) $$\sigma(n_i) \simeq \frac{\text{Ln} 2}{\text{Ln} 3} \text{Log}_2(n_i)$$

The relationship (4.1) can be used to predict the size of $n_i$ required for a given stopping time $\sigma(n_i)$ as the algorithm progresses. Since $n_i$ grows exponentially, it reaches extremely large values ($n_{160}$ in Table 6 has 956 digits). But this does not imply that similar values of $\sigma(n)$ could not be found for integers $n_i$ in lowers regions of the number line. Such a possibility is not precluded by the algorithm which, by design, always searches for greater integers $n_i$.

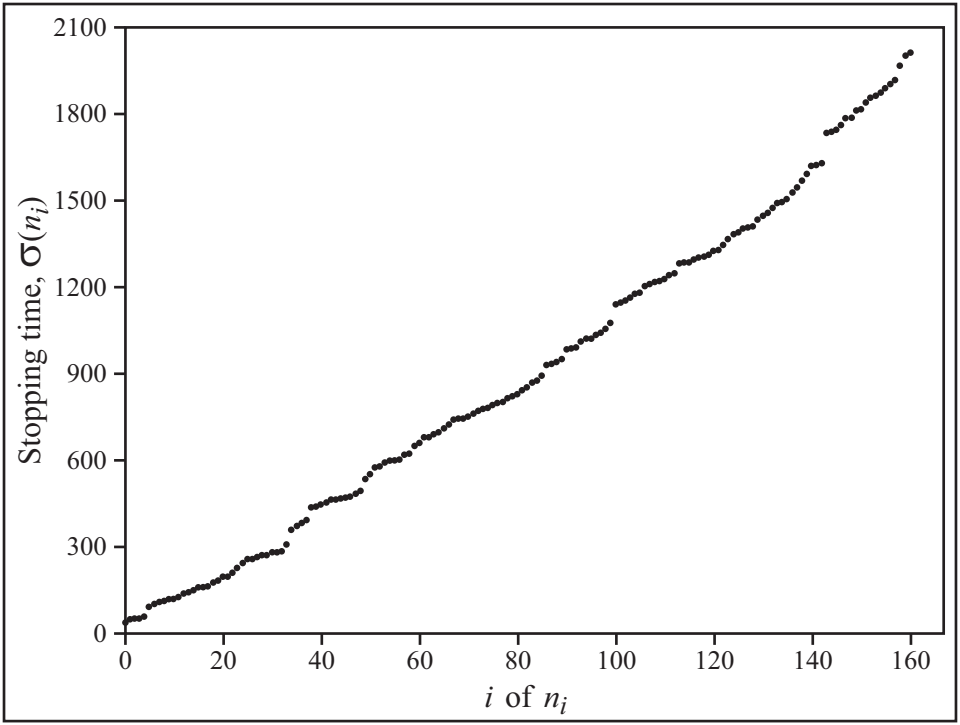

FIGURE 3. Dependence of stopping times $\sigma(n_i)$ on the number of iterations.

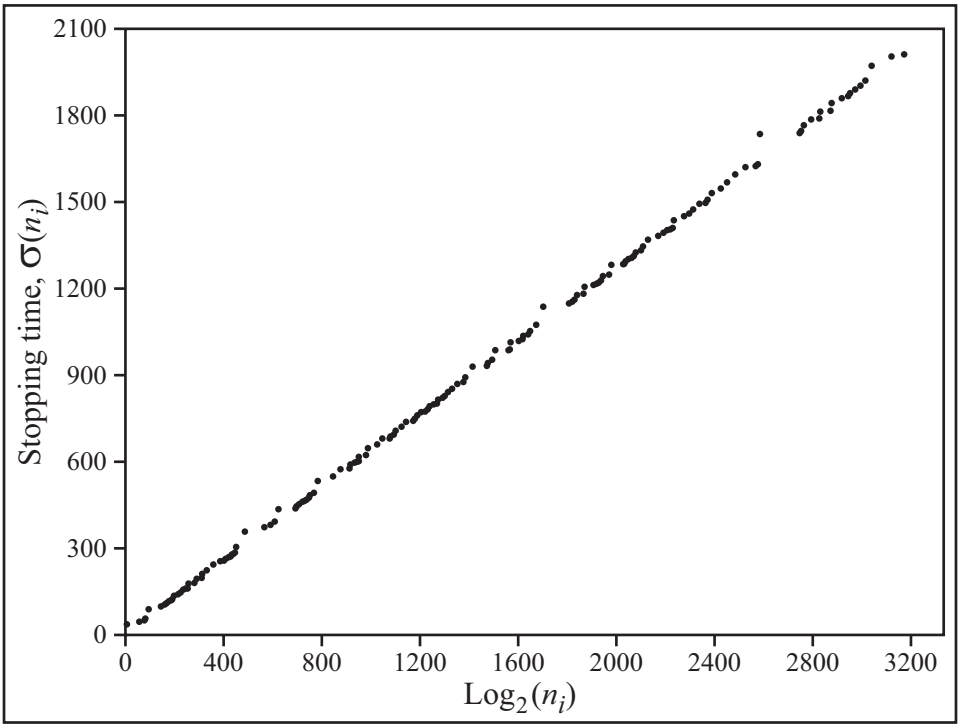

FIGURE 4. Linear relationship between $\sigma(n_i)$ and $\log_2(n_i)$.

## 5. Conclusion

For Collatz Conjecture to hold true, it is required that all positive integers have a finite stopping time. But the results recorded in this article strongly suggest this is not the case. The evidence presented here implies the existence of a infinite set of integers with non-finite stopping times. However, such a conclusion does not contradict preceding results [3,5], since the set of integers with non-finite stopping times will still amount to an infinitesimal fraction of the total set of natural numbers.

To date, no formal proof of Collatz Conjecture has been found. A very simple reason could provide the explanation: the conjecture is false.

## Note

The data presented in Sections 2 and 3 were processed using Excel worksheets. The data included in Section 4 involved very large integers, hence a program written with Mathematica (version 8) was used to generate Table 6.

August 2017
Juan A. Perez. Berkshire, UK.
*Email address*: jap717@juanperezmaths.com